\documentclass[12pt]{article}
\begin{document}
\def\pni{\par \noindent}
\def\vsh{\smallskip}
\def\vs{\medskip}
\def\vvs{\bigskip}
\def\vvvs{\bigskip\medskip} 
\def\vsp{\vsh\pni}
\def\vsn{\vsh\pni}
\def\cen{\centerline}
\def\ra{\item{a)\ }} \def\rb{\item{b)\ }}   \def\rc{\item{c)\ }}
 \cen{{\bf FRACALMO PRE-PRINT} \    {\bf www.fracalmo.org}}
\vsh
\cen{\bf Journal of Computational and Applied Mathematics}
\vsh
\cen{\bf Vol. 153  (2003), pp. 331-342}
\vs
\hrule
   \vskip 0.30truecm

\font\title=cmbx12 scaled\magstep2
\font\bfs=cmbx12 scaled\magstep1
\font\little=cmr10
\begin{center}

{\title Salvatore Pincherle}
\vs

{\title the pioneer of the Mellin-Barnes integrals}

\vvs

 {Francesco MAINARDI}$^{(1)}$, {Gianni PAGNINI}$^{(2)}$

\vs

$\null^{(1)}$
 {\little Department of Physics, University of Bologna, and INFN,} \\
{\little Via Irnerio 46, I-40126 Bologna, Italy} \\
{\little Corresponding Author; E-mail: {\tt francesco.mainardi@unibo.it}}  
\\ [0.20 truecm]
$\null^{(2)}${\little ENEA: Italian Agency for  New Technologies,
  Energy and the Environment} 
\\ {\little Via Martiri di Monte Sole 4, I-40129 Bologna, Italy}\\
{\little E-mail: {\tt gianni.pagnini@bologna.enea.it}}

\end{center} 
\vs

\cen{\bf Abstract} 
\vs
\noindent
The 1888 paper by Salvatore Pincherle (Professor of Mathematics
at the University of Bologna)  on
generalized  hypergeometric functions  is revisited.
We  point out the  pioneering contribution
of the Italian mathematician towards the Mellin-Barnes
integrals based on the duality principle between linear differential
equations and linear difference equation with rational coefficients.
By extending the original arguments used by Pincherle,
we  also show  how to formally derive  the linear differential equation and
the Mellin-Barnes integral representation of the
Meijer $G$ functions.




\newcommand{\FTS}[2]{\frac{{\textstyle #1}}{{\textstyle #2}}}
\def\NN{\hbox{\bf N}}
\def\RR{\hbox{\bf R}}
\def\CC{\hbox{\bf C}}
\def\ZZ{\hbox{\bf Z}}
\font\bfs=cmbx10 scaled \magstep2
\def\eg{{\it e.g.}\ } \def\ie{{\it i.e.}\ }
\def\sg{\hbox{sign}\,}
\def\sgn{\hbox{sign}\,}
\def\sign{\hbox{sign}\,}
\def\e{\hbox{e}}
\def\exp{\hbox{exp}}
\def\ds{\displaystyle}
\def\dis{\displaystyle}
\def\q{\quad}	 \def\qq{\qquad}
\def\lan{\langle}\def\ran{\rangle}
\def\l{\left} \def\r{\right}
\def\lra{\Longleftrightarrow}
\def\arg{\hbox{\rm arg}}
\def\d{\partial}
 \def\dr{\partial r}  \def\dt{\partial t}
\def\dx{\partial x}   \def\dy{\partial y}  \def\dz{\partial z}
\def\rec#1{{1\over{#1}}}
\def\log{\hbox{\rm log}\,}
\def\erf{\hbox{\rm erf}\,}     \def\erfc{\hbox{\rm erfc}\,}


\subsection*{1. Preface}

In Vol. 1, p. 49 of
{\it Higher Transcendental Functions}  of the Bateman Project
\cite{Erdelyi HTF}
  we read
"Of all integrals which contain
gamma functions in their
integrands the most important ones are the so-called Mellin-Barnes
integrals. Such integrals were first introduced by S. Pincherle,
in 1888  \cite{Pincherle 88};
their theory
has been developed in 1910
by H. Mellin (where there are references
to earlier work) \cite{Mellin 10}  and they were used for a complete
integration of the hypergeometric differential equation
by E.W. Barnes \cite{Barnes 08}."
\newpage
\noindent
In the classical treatise on Bessel functions by
Watson \cite{Watson BESSEL}, p. 190, we read
"By using integrals of a type introduced by Pincherle and Mellin,
Barnes has obtained representations of Bessel functions ...."
\vsp
Salvatore Pincherle  (1853 -- 1936)
was Professor of Mathematics at the
University of Bologna from 1880 to 1928.
He retired from the University just after the International
Congress of Mathematicians that he had organized in Bologna,
following the invitation received at the previous Congress
held in Toronto in 1924.
He wrote several treatises and lecture notes on Algebra, Geometry,
Real and Complex Analysis.
His main book related to his scientific activity is entitled
"Le Operazioni Distributive e loro Applicazioni all'Analisi";
it was written in collaboration
with his assistant, Dr. Ugo Amaldi, and
was published in 1901 by Zanichelli, Bologna.
Pincherle can be considered one of the most prominent
founders of the Functional Analysis, as pointed out by J. Hadamard
in his review lecture "Le d\'eveloppement et le
r\^ole scientifique du Calcul fonctionnel",
given at the Congress of Bologna (1928).
A description of Pincherle's scientific works requested from him
by Mittag-Leffler, who was the Editor of Acta Mathematica,
appeared (in French) in 1925 on  this prestigious
journal \cite{Pincherle TRAVAUX}.
A collection of selected papers (38 from 247 notes plus 24
treatises) was edited by Unione Matematica Italiana (UMI)
on the occasion of the centenary of his birth,
and published by Cremonese, Roma 1954.
Note that S. Pincherle was the first President of UMI, from 1922 to 1936.
Here we point out that the 1888 paper (in Italian)
of S. Pincherle  on the
{\it Generalized Hypergeometric Functions}
led him to  introduce the afterwards named Mellin-Barnes integral
to represent the solution of a
generalized hypergeometric differential equation  investigated
by Goursat  in 1883.
Pincherle's priority was
explicitly recognized by Mellin and Barnes themselves,
as reported below.
\vsp
In 1907 Barnes, see p. 63 in \cite{Barnes 07},	wrote: 
"The idea of employing contour integrals involving gamma functions
of the variable in the subject of integration appears to be due
to Pincherle, whose suggestive paper was the starting point of
the investigations of Mellin (1895)  though the type of contour
and its use can be traced back to Riemann."
In 1910 Mellin, see p. 326ff in  \cite{Mellin 10},
devoted a section (\S 10: Proof of Theorems of Pincherle)
to revisit the original work of Pincherle; in particular,  he wrote
"Before we are going to prove this theorem, which is a special case of a
 more general theorem of Mr. Pincherle, we want to describe more closely
 the lines $L$ over which the integration preferably is to be carried
 out." [free translation from German].
\newpage
\noindent
The Mellin-Barnes integrals are the essential tools
for treating the
two classes of higher transcendental functions known as $G$ and $H$
functions, introduced by Meijer (1946) \cite{Meijer 46}
and Fox (1961) \cite{Fox 61} respectively, so
Pincherle can be considered their precursor.
For an exhaustive treatment of the Mellin-Barnes integrals we refer
to the recent monograph by Paris and Kaminski
\cite{ParisKaminski BOOK01}.
\vsp
The purpose of our paper is to let know  the community
of scientists interested in special functions the
pioneering 1888 work by Pincherle, that, in the author's
intention, was devoted
to compare two different generalizations of the
Gauss hypergeometric function due to Pochhammer and to
Goursat. Incidentally, for a particular case of the Goursat function,
Pincherle used an integral representation in the complex plane
that in future was adopted by Mellin and Barnes for their
treatment     of the {\it generalized hypergeometric functions}
known as $\,_pF_q (z)\,. $
We also intend to show,
in the original part of our paper,
that, by extending the original arguments
by Pincherle, we are able to provide
the Mellin-Barnes integral representation of the transcendental
functions introduced by Meijer (the so-called $G$ functions).
\vsp
The paper is divided as follows.
In Section 2
 we report the major statements and
results of the	1888 paper by Pincherle.
In Section 3 we show how it is possible to originate from these
results  the Meijer $G$ functions by a proper 
generalization of  Pincherle's method.
Finally, Section 4 is devoted to the conclusions.
We find it convenient to reserve an Appendix for recalling
some basic notions for the generalized hypergeometric functions
and the Meijer $G$ functions.


\subsection*{2. The Pochhammer and Goursat generalized hypergeometric
functions via Pincherle's  arguments}

The 1888 paper by Pincherle is based on what he called
"duality principle",
which relates linear differential equations with
rational coefficients to linear difference equations with
rational coefficients.
Let us remind that the sentence "rational coefficients" means
that the coefficients are in general rational functions
(\ie ratio between two polynomials)
of the independent variable and, in particular, polynomials.
\vsp
By using  this principle with polynomial coefficients,
Pincherle  showed that
two generalized  hypergeometric functions
proposed and investigated respectively by   Pochhammer (1870), see
\cite{Pochhammer 70}, and by Goursat (1883),
see \cite{Goursat CR83,Goursat ENS83},
can be obtained and related to each other
\footnote{In fact, translating from Italian,
 the author so writes in introducing his paper:
 "It is known that to any linear differential
 equation with rational coefficients one may let correspond
 a linear difference equation with rational coefficients.
 In other words, if the former equation is given, one can immediately
 write the latter one and viceversa;
 furthermore, from the integral of the one,
 the integral of the latter can be easily deduced.
 This {\it relationship} appears to be originated by a sort of
 {\it duality principle} of which, in this note, I  want to treat
 an application concerning generalized hypergeometric functions."}.
\vsp
  The generalized hypergeometric functions introduced by Pochhammer
and Goursat   considered by Pincherle are solutions
of linear differential equations of order $n$
with polynomial coefficients, that we report in Appendix.
As a matter of fact, the duality principle states the correspondence
between a  linear {\it ordinary differential equation} ($ODE$)
and a linear {\it finite difference equation} ($FDE$).
The  coefficients of both  equations are assumed to be
{\it rational} functions, in particular {\it polynomials}.
In his analysis \cite{Pincherle 88} Pincherle considered
the correspondence between the following equations,
$$
  \sum_{h=0}^m
\l( a_{h\,0} + a_{h\,1}\,\e^{-t}
  + a_{h\,2}\,\e^{-2t}	+ \dots + + a_{h\,p}\,\e^{-pt}	  \r)
 \psi^{(h)} (t) =0, \eqno(2.1)$$
$$   \sum_{k=0}^p
\l[ a_{0\,k} + a_{1 \,k}   (x +k)
  + a_{2\,k}  (x +k)^2	+ \dots +  a_{m\,k}   (x +k)^m	  \r]
 f (x +k) =0,\eqno(2.2)$$
where $\psi(t)$ and $f(x)$ are analytic functions.
  These functions are required to be related to each other
 through a	Laplace-type transformation
$ \psi(t) \,{\leftrightarrow} \, f(x)$ defined by the formulas
$$ \hbox{(a)} \q
     f(x) = \int_l \e^{-x t}\, \psi (t)\, dt\,,
\qq
\hbox{(b)} \q
 \psi (t)=\rec{2\pi i}\,\int_L \e^{+x t}\, f (x)\,dx \,,
\eqno(2.3)$$
where $l$ and $L $ are	appropriate integration paths  in the
complex $t$ and $x$ plane, respectively.
\vsp
The singular points of the $ODE$ are the roots of the polynomial
$$
  a_{m\,0} + a_{m\,1}\,z
  + a_{m\,2}\,z^2  + \dots +  a_{m\,p}\,z^p =0\,.
  \eqno(2.4)$$
whereas the singular points of the $FDE$ are the rots of
the polynomial
$$
 a_{0\,0} + a_{1 \,0}	z
  + a_{2\,0}  z^2  + \dots +  a_{m\,0}	 z^m = 0\,.
 \eqno(2.5)$$
For  the details of the above correspondence
Pincherle refers to the 1885 fundamental paper by Poincar\'e
 \footnote{For an account of Poincar\'e's theorem upon which
Pincherle based his analysis the interested reader can consult
the recent book  by Elaydi  \cite{Elaydi DE99}, pp. 320-323.}
 \cite{Poincare 85},
 and his own 1886 note \cite{Pincherle 86}.
Here we limit ourselves to point out what can be easily seen
from a formal comparison between the $ODE$ (2.1) and the $FDE$ (2.2).
We recognize
that the degree $p$ of the coefficients in $\e^{-t}$ of the $ODE$
provides the order of the $FDE$,  and that the order $m$
of the $ODE$  gives the degree in $x$
of the coefficients of the $FDE$.
Viceversa, the degree $m$ of the coefficients  of the $FDE$ provides
the order of the $ODE$, and the order $p$ of the $FDE$
gives the degree in $\e^{-t}$ of the coefficients
of the $ODE$.
\vsp
Pincherle's intention was to apply the above duality principle
in order to
compare the generalized
hypergeometric function introduced by Pochhammer and governed by
(A.7) with that by  Goursat governed by (A.6).
Using his words, he proved that the family of the Pochhammer functions
 (of arbitrary order $p$)  originates
from a linear $FDE$ (of  order $p$) whose coefficients are polynomials
of the first degree in $x\,,$
and that the family of the Goursat functions (of arbitrary order $m$)
originates from a linear $ODE$ (of order $m$) whose
coefficients are polynomials of the first degree in $x = \e^{-t}\,. $
As a consequence of the duality principle there is a
mutual correspondence between the properties of the functions
belonging to one family and to the other.
\vsp
For the Pochhammer function he started from
the $ODE$ of the first order
   $$ \qq \qq \l( a_{0\,0} + a_{0 \,1}\, \e^{-t}
 + a_{0 \,2}\, \e^{-2 t}
  + \dots  + a_{0\,p}\, \e ^{-pt}    \r)  \,\psi(t) \qq\qq\qq\qq
						    \eqno(2.6)$$
$$ +  \l( a_{1\,0} + a_{1 \,1}\, \e^{-t}
+ a_{1 \,2}\, \e^{-2 t}
  + \dots  + a_{1\,p}\, \e ^{-pt} \r)  \,\psi^{(1)}(t)
    = 0\,, $$
to be put in correspondence with the $FDE$
$$ \l( a_{0\,0} + a_{1 \,0} x\r) f(x)+
  \l[ a_{0\,1} + a_{1 \,1} (x+1)\r] f(x+1)+
  \l[ a_{0\,2} + a_{1 \,2} (x+2)\r] f(x+2)  \eqno(2.7)$$
$$ +  \dots  +	 \l[ a_{0\,p} + a_{1 \,p} (x+p)\r] f(x+p) =0\,.
$$
In this case Pincherle was able to show that the solution $f(x)$
of the $FDE$ (2.7), obtained through the  formula (a) in (2.3),
depends on $p$ parameters, whose logarithms are the singular
points of the $ODE$ (2.6). With respect to each of these
parameters, $f(x)$ satisfies a linear $ODE$ of the Pochhammer type
of order $p\,.$ 
\vsp
For the Goursat function he started from a $FDE$ of the first order
$$ \qq \qq \qq \qq \l[ a_{0\,0} + a_{1 \,0}\, x    +  a_{2 \,0}\, x^2
  +   \dots  + a_{m\,0}\, x ^m	  \r]  \,f(x )\qq\qq\qq
 \eqno(2.8)$$
$$ + \l[ a_{0\,1} + a_{1 \,1}\,(x+1)  + a_{2 \,1}\,(x+1)^2
  +   \dots  + a_{m\,1}\, (x+1) ^m    \r]  \,f(x+1) =0\,,$$
to be put in correspondence to the linear $ODE$ of order $m$
$$\qq  \l( a_{0\,0} + a_{0 \,1} \e^{-t}\r)\, \psi(t) +
 \l( a_{1\,0} + a_{1 \,1} \e^{-t}\r)\, \psi^{(1)}(t) +
  \l( a_{2\,0} + a_{2 \,1} \e^{-t}\r)\, \psi^{(2)}(t)
\qq \eqno(2.9)$$
$$ + \dots
+ \l( a_{m\,0} + a_{m \,1} \e^{-t}\r)\, \psi^{(m)}(t) =0\,.
$$
Using a result of  Mellin, see \cite{Mellin 86,Mellin 87},
Pincherle wrote
the solution of the $FDE$ (2.8) as
$$ f(x ) = c^x	 \prod_{\nu =1}^m
     {\Gamma(x	-\rho _\nu )\over \Gamma(x -\sigma _\nu)}\,,
\eqno(2.10)$$
where
the $\rho _\nu$'s and $\sigma _\nu$'s are respectively the roots of the
algebraic  equations
$$\cases{ a_{0\,0} + a_{1 \,0}\, x
   + \dots  + a_{m\,0}\, x ^m  = a_{m\,0}\,
{\ds \prod_{\nu =1}^m}
     (x  -\rho _\nu )= 0\,,\cr\cr
  a_{0\,1} + a_{1 \,1}\, (x +1)
  + \dots +  a_{m\,1}\, (x +1)^m =
 a_{m\,1}\, {\ds \prod_{\nu =1}^m}
     (x  -\sigma  _\nu ) =0\,.
   \cr}
\eqno(2.11)
$$
and $c$ is a constant. If $a_{m\,0},\,a_{m\,1}$ are both
different from zero, we can assume $c = - a_{m\,0}/ a_{m\,1}\,.$
\vsp
Pincherle showed that, by setting $z= c\,\e^{t}\,,$
the $ODE$ of order $m$ (2.9) is nothing but
the Goursat differential equation (A.6).
\vsp
Furthermore, in the special case $a_{m\,1} =0\,,$ he gave the following
relevant formula for the solution
$$ \psi(t) = \rec{2\pi i} \,\int_{a-i\infty} ^{a+i\infty}
 {\Gamma(x-\rho_1)\,\Gamma(x-\rho_2)\dots \Gamma(x-\rho_m)
 \over
 \Gamma(x-\sigma _1)\,\Gamma(x-\sigma_2)\dots \Gamma(x-\sigma_{m-1})}
\, \e^{xt} \, dx \eqno(2.12)$$
where $a > \Re \{\rho _1, \rho _2, \dots, \rho _m\}\,. $
We recognize in (2.12) the first example in the literature
of the (afterwards named)  Mellin-Barnes integral.
\vsp
The convergence of the integral was proved by Pincherle
by using his asymptotic formula for $\Gamma(a + i\eta)$
as $\eta \to \pm \infty$
\footnote{We also note the priority of Pincherle in
obtaining this asymptotic formula, as outlined by Mellin, see \eg
\cite{Mellin 91}, pp. 330-331, and  \cite{Mellin 10}, p.309. 
In his 1925 "Notices sur les travaux"
\cite{Pincherle TRAVAUX} (p. 56, \S 16) Pincherle wrote
"Une expression asymptotique de $\Gamma(x)$ pour $x \to \infty$
dans le sens imaginaire qui se trouve dans
\cite{Pincherle 88} a \'et\'e attribu\'ee \`a
d'autres auteurs, mais M. Mellin m'en a
r\'ecemment r\'evendiqu\'e  la priorit\'e."
This formula is fundamental to investigate the convergence of the
Mellin-Barnes integrals, as one can recognize from the detailed
analysis by Dixon and Ferrar \cite{DixonFerrar 36}, see also
\cite{ParisKaminski BOOK01}.}.
So,  for a solution of a particular case of the Goursat equation,
Pincherle provided an integral representation
that later was adopted by Mellin and Barnes for their
treatment   of the generalized hypergeometric functions
$\,_p F_q (z)\,.$
Since then, the merits of Mellin and Barnes were so well recognized
that their names were attached to the integrals
of this type; on the other hand,
after the 1888 paper (written in Italian),
Pincherle did not pursue on this topic,
so his name was no longer related
to   these integrals and, in this respect, his 1888 paper
was practically ignored.

\subsection*{3. The Meijer transcendental function via Pincherle's arguments}

In more recent times
other families of higher transcendental functions
have been introduced to generalize the hypergeometric function
based on their representation by Mellin-Barnes type integrals.
We especially refer to the so-called  $G$ and $H$  functions,
briefly recalled in the Appendix.
\vsp
In this section (the original part of our paper)
we show that by extending the original arguments  by Pincherle
based on the duality principle	we are able to provide
the differential equation and the
Mellin-Barnes integral representation of the $G$ functions.
However, we note  that these arguments, being
based on equations with rational coefficients,
do not allow  us to treat the Fox $H$ functions,
since for them an ordinary differential equation
cannot be found in the general case.
\vsp
Our starting point is still the "duality principle"
that involves a {\it $FDE$ of the  first order
as in Pincherle's approach for the Goursat function},
but, at variance of Eq. (2.8), we now allow
that the degree of the two 
polynomial coefficients are not necessarily equal.
Setting $p, q$ the degrees of these coefficients,  our $FDE$
reads
$$\qq \qq \qq \qq \l[ a_{0\,0} + a_{1 \,0}\, x	  +  a_{2 \,0}\, x^2
  +   \dots  + a_{p\,0}\, x ^p \r]  \,f(x )\qq \qq \qq\eqno(3.1)$$
$$ + \l[ a_{0\,1} + a_{1 \,1}\,(x+1)  + a_{2 \,1}\,(x+1)^2
  +   \dots  + a_{q\,1}\, (x+1) ^q    \r]  \,f(x+1) =0\,.$$
We can prove after some algebra
that the associated $ODE$
turns out to be  independent of the
order relation between $p$ and $q$ and reads
$$
  \sum_{h=0}^p	 a_{h\,0} \, \psi^{(h)} (t)  +
    \e^{-t} \,	 \sum_{h=0}^q  a_{h\,1}\,
 \psi^{(h)} (t) =0\,.\eqno(3.2)$$
As we have learnt from Pincherle's analysis,
the solution $\psi(t)$	of the $ODE$ (3.2) can
be expressed  in terms
of the solution $f(x)$ of the $FDE$ (3.1),  according to
the integral representation  (b) in Eq. (2.3).
\vsp
Now, in view of Mellin's  results  used by Pincherle
(see also Milne-Thomson \cite{MilneThomson FDE51}, \S 11.2, p. 327),
we can write
the solution of (3.1) in terms of products and fractions
of $\Gamma$ functions.
Denoting  by
$\rho _j$ ($j=0,1,\dots ,p$) and $\sigma_k$ ($k=0,1,\dots,q$)
the roots of the
algebraic  equations
$$\cases{ a_{0\,0} + a_{1 \,0}\, x
   + \dots  + a_{p\,0}\, x ^p  = a_{p\,0}\,
{\ds \prod_{j =1}^p}
     (x  -\rho _j )= 0\,,\cr\cr
  a_{0\,1} + a_{1 \,1}\, (x +1)
  + \dots +  a_{q\,1}\, (x +1)^q =
 a_{q\,1}\, {\ds \prod_{k =1}^q}
     (x  -\sigma_k ) =0\,.
   \cr}
\eqno(3.3) $$
we can write the required solution as
$$ f(x) = c^x\,
    {\prod_{j =1}^p \Gamma(x-\rho_j)\over
     \prod_{k =1}^q \Gamma(x -\sigma_k)   }\,,
 \q c = - { a_{p\,0} \over a_{q\,1}}\,.
  \eqno(3.4)   $$
We note, by using the known properties
of the Gamma function,
that Eq. (3.4) can be re-written
in the following alternative form
$$ f(x) = c^x\,
    {\prod_{k =1}^q \Gamma(1+\sigma_k -x)\over
     \prod_{j =1}^p \Gamma(1 +\rho_j -x)      }\,,
 \q c =  (-1)^{p-q+1}\, { a_{p\,0} \over a_{q\,1}}\,.
  \eqno(3.5)   $$
Furthermore,
introducing the integers $m,n$ such that
$ \,0\le m \le q\,,$ $\, 0\le n \le p\,, $
we can combine the previous formulas (3.4)-(3.5)
and obtain the alternative form 
$$ f(x) = c^x\,
    {\prod_{j =1}^n \Gamma(x-\rho_j)\,
     \prod_{k =1}^m \Gamma(1+\sigma_k -x)
	    \over
     \prod_{j =n+1}^p \Gamma(1+\rho_j -x)\,
     \prod_{k =m+1}^q \Gamma(x- \sigma_k)}  \,,  \eqno(3.6)$$
with
$$ c =	(-1)^{m+n -p+1}\, { a_{p\,0} \over a_{q\,1}}\,.
  \eqno(3.7)   $$
We note that Eq.  (3.6) reduces to the Pincherle expression (2.10)
by setting $\{n=p=q\,,\, m=0\}$, and  to Eqs (3.4), (3.5)
by setting $\{n=p\,,\,m=0\}$,
$\{n=0\,,\,m=q\}$,  respectively.
By adopting the form (3.6)-(3.7) we have the most general
expression for $f(x)$ which in its turn allows us to arrive
at the most general solution $\psi(t)$
of the corresponding $ODE$  (3.2) in the form
$$ \psi(t)
 =\rec{2\pi i}\,\int_L c^x  \,
 {\prod_{j =1}^n \Gamma(x-\rho_j)\,
     \prod_{k =1}^m \Gamma(1+\sigma_k -x)
	    \over
     \prod_{j =n+1}^p \Gamma(1+\rho_j -x)\,
     \prod_{k =m+1}^q \Gamma(x- \sigma_k)}  \,
  \e^{xt}\, dx \,,
\eqno(3.8)$$
where  $L$ is an  appropriate integration path	in the
complex  $x$ plane.
\vsp
Now, starting from  (3.2)  and (3.7)-(3.8) it is
not difficult to arrive at the general $G$ function
namely at  its $ODE$ 
and at its Mellin-Barnes integral representation,
both given in Appendix. 
For this purpose we need only to carry out some algebraic manipulations
and obvious transformations of variables.
\vsp
We first note that using (3.3) the $ODE$ (3.2) reads
$$ \l[ a_{p\,0}\, \prod_{j=1}^p
\l( {d\over dt} -\rho _j\r) +
 a_{q\,1}\,\e^{-t}\,  \prod_{k=1}^q
\l({d\over dt} -\sigma_k -1\r)\r]\,\psi(t)=0\,. \eqno(3.9)$$
Then, putting
$$ z = c\, \e^t\,, \q u(z) = \psi[t(z)]\,,
\q a_j= 1 + \rho _j\,,\q b_k= 1+ \sigma_k\,, \eqno(3.10)$$
and using (3.7), we get  from (3.9)
$$ \l[(-1)^{p-m-n}z
\prod_{j=1}^p\l(z{d \over dz}-a_j+1\r)-
\prod_{k=1}^q\l(z{d \over dz}-b_k\r)\r]
 u(z)=0 \,, \eqno(3.11)$$
which is just the $ODE$ satisfied by the Meijer $G$ function
of orders $m,n,p,q$, see (A.10).
Of course, at least formally, the Mellin-Barnes
integral
representation of the $G$ function (A.8)-(A.9)
is recovered as well and reads
(setting $s = x$)
$$ u(z) = \rec{2\pi i}\,
\int_L
{
 \prod_{k=1}^m \Gamma(b_k- s)\, \prod_{j=1}^n \Gamma(1-a_j + s)
\over
\prod_{k=m+1}^q \Gamma(1-b_k + s)\,\prod_{j=n+1}^p \Gamma(a_j- s)
}
\, z^s \, ds\,.\eqno(3.12)$$

\newpage
\subsection*{4. Conclusions}

We have revisited the 1888 paper (in Italian) by Pincherle
on generalized	hypergeometric functions,
based on the duality principle between linear differential
equations and linear difference equation with rational coefficients.
We have  pointed out the  pioneering contribution
of the Italian mathematician towards the afterwards named Mellin-Barnes
integral representation that he was able to provide
for a special case of a generalized hypergeometric
function introduced by Goursat in 1883.
By extending his original arguments
we  have shown	how to formally derive	the ordinary differential
equation and
the Mellin-Barnes integral representation of the
$G$ functions introduced by Meijer in 1936-1946.
 So, in principle,  Pincherle could have introduced
the $G$ functions much before Meijer if he had intended to
 pursue his original arguments in this direction.
Finally, we like  to point out that
the so-called Mellin-Barnes integrals are an efficient
tool to deal with  the higher transcendental functions.
In fact,  for a pure mathematics view point
they facilitate the representation 
of these functions
(as formerly indicated by Pincherle),
and for an applied mathematics	view point
they can be successfully adopted to compute
the same functions.
In this respect we refer to the recent paper by
Mainardi, Luchko and Pagnini \cite{Mainardi LUMAPA01},
who have computed the solutions of diffusion-wave equations
of fractional order by using their Mellin-Barnes
integral representation.

\subsection*{Acknowledgements}
 \noindent
Research performed under the auspices of the National Group of Mathematical
Physics (G.N.F.M. - I.N.D.A.M.) and partially supported
 by the Italian Ministry
of University (M.I.U.R) through the Research Commission of the
 University of Bologna and by the National Institute of Nuclear
Physics (INFN)	through the Bologna branch (Theoretical Group).
The authors are grateful to Prof. R. Gorenflo for the  discussions
and  the helpful comments.

\subsection*{Appendix: Some generalizations of the hypergeometric functions}

The purpose of this Appendix is to provide a survey of
some higher transcendental functions that have been proposed
for generalizing the 
hypergeometric function.
In particular we shall consider the functions  investigated by Pochhammer
(1870) and Goursat (1883), that have interested Pincherle
in his 1888 paper, and	the $G$ functions introduced by Meijer
(1936-1946), since they are re-derived in our
present analysis  by extending the 
arguments by Pincherle.
Our survey is essentially based on the classical  handbook
of the Bateman Project
\cite{Erdelyi HTF} and on the more
recent treatise by Kiryakova \cite{Kiryakova 94}.
\vsp
Let us start by recalling  the classical  {\it hypergeometric
equation}.
If a homogeneous linear differential equation  of the second order
has at most three singular points we may assume that these
are $0, 1, \infty\,. $ If all these singular points are "regular", then
the equation can be reduced to the form
$$ z(1-z) \,   {d^2 u  \over dz^2} +
  [c-(a +b +1) z ] \, {d u  \over dz} - ab \, u(z)
=0\,,
\eqno(A.1)
$$
where $a,b,c$ are arbitrary complex constants.
This is the {\it hypergeometric equation}.
Taking	$c \ne 0,-1,-2,\dots\,,$
and  defining the Pochhammer symbol
$$ (\alpha )_n = {\Gamma(\alpha +n)\over \Gamma(\alpha )}\,,
\; \hbox{\ie} \;
(\alpha )_0=1,,\;
 (\alpha)_n =\alpha (\alpha+1)\dots (\alpha +n-1)\,,\; n=1,2,\dots$$
 then	the solution of Eq. (A.1), regular at $z=0\,, $
known as {\it Gauss hypergeometric function},
turns out to be
$$ u(z) = \sum_{n=0}^\infty
{ (a)_n\, (b)_n\over (c)_n n!}\, {\ds z^n} := F(a,b;c;z)
\,.\eqno(A.2)$$
\vsp
The above hypergeometric series can be generalized by introducing
$p$ parameters $a_1,\dots a_p$ 
(the numerator-parameters) and
 $q$ parameters $b_1, \dots, b_q$ 
(the denominator-parameters).
The ensuing series
$$  \sum_{n=0}^\infty
{ (a_1)_n\,\dots (a_p)_n\over (b_1)_n \dots (b_q)_n}
  \, {\ds {z^n \over n!}}
 \, :=
\; _pF_q (a_1, \dots, a_p ; b_1, \dots, b_q ;z) \, ,\eqno(A.3)$$
 or, in a more compact form,
$$\sum_{n=0}^{\infty}
{\Pi_{j=1}^p (a_j)_n \over \Pi_{k=1}^q (b_j)_n} \,
{\ds {z^n \over n!}} \,:=
\; _pF_q \l[ (a_j)_1^p ; (b_k)_1^q ;z) \r]
  \eqno(A.3')$$
is known as the {\it generalized hypergeometric series}.
In general
  (excepting certain integer values of the parameters
for which the series fails to make sense or terminates
\footnote{If at least one of the denominator parameters $b_k$
($k =1,\dots, q$) is zero or a negative integer, Eq. (A.3) has no meaning
at all, since the denominator of the general term vanishes
for a sufficiently large index. If some of the numerator
parameters are zero or negative integers, then the series terminates
and turns into a {\it hypergeometric polynomial}.})
the series $\, _pF_q$	 converges for all finite $z$ if $p\le q\,,$
converges for $|z|<1$ if $p=q+1\,,$ and diverges for
all $z \ne 0$ if $p>q+1\,. $
The resulting {\it generalized hypergeometric function}
$u(z) =\, _pF_q$
will
satisfy a   {\it generalized hypergeometric equation}.
If we note that Eq. (A.1) satisfied by $u(z) =	F(a,b;c;z) =
\,_2F_1(a,b;c;z)   $
can be written in the equivalent form
(see \eg Rainville \cite{Rainville SF60}, \S 46, p. 75) :
$$ \l[z
\l(z{d \over dz}+a\r)\l(z{d \over dz} +b \r)-
z{d \over dz} \l(z{d \over dz}+c-1\r)\r]
u(z)=0 \,, \eqno(A.1')$$
we arrive at the  equation of order $n=q+1$ for $u(z)
=\, _pF_q \l[ (a_j)_1^p ; (b_k)_1^q ;z) \r]\,:$
$$ \l[z
\prod_{j=1}^p\l(z{d \over dz}+a_j\r)-
z{d \over dz}\prod_{k=1}^{q}\l(z{d \over dz}+b_k-1\r)\r]
u(z)=0 \,. \eqno(A.4)$$
The above equation containing the operator $z d/dz$ can be
written in a more explicit form by using $D = d/dz$,
see \eg \cite{Erdelyi HTF} \S 42, p.184.
Distinguishing between the cases $p\le q$
and $p=q+1\,, $ we get
 the following general equations
in $v= v(z)\,:$
$$ z^q D^{q+1} v +
\sum_{\nu =1}^q z^{\nu -1} (A_\nu  z - B_\nu) \, D^\nu v + A_0 v = 0\,,
 \q p \le q\,, \eqno(A.5)$$
$$  z^q (1-z)D^{q+1} v +
\sum_{\nu =1}^q z^{\nu -1} (A_\nu  z - B_\nu ) \, D^\nu  v
 + A_0 v  = 0\,,
  \q p =q+1\,, \eqno(A.6)$$
where $A_0, A_\nu, B_\nu$ are constants.  
Eq. (A.5) has two singular points, $z=0, \infty$ of which
$z=0$ is of regular type, whereas
Eq. (A.6) has	three singular points, $z=0,1, \infty$ of
regular type, like Eq. (A.1).
An equation of the same type as  Eq. (A.6) was formerly introduced in
1883 by Goursat \cite{Goursat CR83,Goursat ENS83}
in his essay on hypergeometric
functions of higher order.
\vsp
Another generalization of the Gauss hypergeometric  equation
was previously proposed in 1870
by Pochhammer \cite{Pochhammer 70}.
He 
investigated  the most general homogeneous linear differential equation
of the order $n$ ($n>2$) of Fuchsian type, namely with only "regular"
singular points in  $\{a_1, a_2, \dots, a_n,\infty\}\,. $
The Pochhammer function thus satisfies a differential equation	of
the type
$$  \phi _n(z) \, {d^n w  \over dz^n} +
  \dots  + \phi _1(z) \, {d w  \over dz} +    \phi_0\, w(z) = 0\,,
\eqno(A.7)
$$
where the coefficients $\phi _\nu (z) $ ($\nu  = 0, 1 \dots, n$) are
polynomials of degree $\nu \,, $ with
$\phi _n(z) = (z-a_1)(z-a_2)\dots (z-a_n)\,. $
\vsp
The $\,_pF_q\,$ functions  satisfying  Eqs (A.5)-(A.6)
and the Pochhammer functions satisfying  Eq. (A.7)
are not the only generalizations of the Gauss hypergeometric
function (A.2). In 1936  Meijer \cite{Meijer 36}
  introduced a new class of transcendental
functions, the so called $G$ functions, which provide
an interpretation of the  symbol  $\,_pF_q$ when
$ p >q+1\,. $  Originally, the $G$ function was defined in a manner
resembling (A.2). Later \cite{Meijer 46}, this definition was replaced
by one in terms of Mellin-Barnes type integrals.
The latter definition has the advantage that it allows
a greater freedom in the relative values of $p$ and $q$.
Here, following \cite{Erdelyi HTF}, we shall complete
Meijer's definition so as to include all values of $p$ and $q$
without placing any (non-trivial) restriction on $m$ and $n$.
One defines
$$   G^{m,n}_{p,q}
\left[ z \left\vert
{a_1, \dots , a_p\atop b_1, \dots , b_q}
 \right.
\right]=
G^{m,n}_{p,q}
\left[ z \left\vert
(a_j)_1^p \atop (b_j)_1^q
 \right.
\right] =
\rec{2\pi i}\,
\int_L {\cal{G}}^{m,n}_{p,q} (s) \, z^s \, ds\,, \eqno(A.8)$$
where
 $ L$ is a suitably chosen path, $z \neq 0\,,$
 $ z^s :=\exp \l[ s(\ln |z| + i \, \hbox{arg} \,z)\r]$
with a single valued branch of $\hbox{arg}\, z$,
and  the integrand is defined as follows
$$ {\cal{G}}^{m,n}_{p,q} (s)=
{
 \prod_{k=1}^m \Gamma(b_k- s)\, \prod_{j=1}^n \Gamma(1-a_j + s)
\over
\prod_{k=m+1}^q \Gamma(1-b_k + s)\,\prod_{j=n+1}^p \Gamma(a_j- s)
}\,.\eqno(A.9) $$
In (A.9) an empty product is interpreted as 1, the integers
$m,n,p,q$ (known as {\it orders of the $G$ function}) are such that
$0\le m \le q\,,$
$\,0 \le n \le p\,,$ and the parameters
$a_j$ and $b_k$ are such that
no pole of $\Gamma(b_k-s), $   $k=1,\dots, m,$
coincides with any pole  of $\Gamma(1-a_j+s), $   $j=1,\dots, n.$
For the details of the integration path,
which can be of three different types,
we refer to \cite{Erdelyi HTF} (see also \cite{Kiryakova 94}
where an illustration of what these contours can be like
is found).
\vsp
One can establish that the Meijer $G$ function $u(z)$
satisfies the linear ordinary differential equation of
generalized hypergeometric type,
see \eg \cite{Kiryakova 94}  [p. 316, Eq. (A.19)],
$$ \l[(-1)^{p-m-n}z
\prod_{j=1}^p\l(z{d \over dz}-a_j+1\r)-
\prod_{k=1}^q\l(z{d \over dz}-b_k\r)\r]
u(z)=0 \,. \eqno(A.10)$$
For  more details on the Meijer function and on the
singular points of the above  differential equation we
 refer to  \cite{Kiryakova 94}. Here we limit ourselves to
show how
the generalized hypergeometric function $\,_pF_q$
can be expressed in terms of a Meijer $G$ function
and  thus in terms of  Mellin-Barnes integral. We have
$$ _pF_q  ( (a)_p ;(b)_q ;z ) =
{\Pi_{k=1}^q \Gamma(b_k) \over
\Pi_{j=1}^p \Gamma(a_j)}
G^{1,p}_{p,q+1}
\left[- z   \left\vert
{\hfill (1-a_j)_{1}^{p} \hfill \atop
 \hfill 0,\, (1-b_k)_{1}^{q} \hfill}
 \right. \right] \,, \eqno(A.11)$$
$$ G^{1,p}_{p,q+1} =
 {1 \over 2 \pi i} \,\int_{- i \infty}^{+ i\infty}
{\Gamma(a_1+s) \cdots \Gamma(a_p+s)\Gamma(-s) \over
\Gamma(b_1+s) \cdots \Gamma(b_q+s)}
(-z)^s \, ds \,, \eqno(A.12) $$
$$ a_j \ne 0,-1,-2,\dots; \q j=1,\dots,p; \q
|\hbox{arg}(1-zi)|<\pi \,. \eqno(A.13)$$
Here the path of integration is the imaginary axis (in the complex
$s$-plane) which can be deformed, if necessary, in order to
separate the poles of $\Gamma(a_j+s)$, $j=1,\dots,p$
from those of $\Gamma(-s)\,. $
\vsp
Though the $G$ functions are quite general in nature, there still
exist examples of special functions,
like the Mittag-Leffler and the  Wright functions,
which do not form their particular cases.
A more general class  which includes
those functions can be achieved by introducing the
Fox $H$ functions \cite{Fox 61}, whose representation in terms of the
Mellin-Barnes integral is a straightforward generalization of that
for the $G$ functions.
For this purpose we need to add to the sets of the  complex
parameters $a_j$ and $b_k$ the new sets of positive numbers $\alpha _j$
and $\beta _k$ with $j=1,\dots, p,$ $k=1,\dots,q,$
and modify in the integral of (A.8) the  kernel
$ {\cal{G}}^{m,n}_{p,q} (s)$  into
$$		 {\cal{H}}^{m,n}_{p,q} (s)=
{
 \prod_{k=1}^m \Gamma(b_k- \beta_k s)\,
 \prod_{j=1}^n \Gamma(1-a_j + \alpha _j s)
\over
\prod_{k=m+1}^q \Gamma(1-b_k + \beta_k	s)\,
\prod_{j=n+1}^p \Gamma(a_j- \alpha_j  s)
}\,.\eqno(A.14) $$
Then the Fox $H$ function turns out to be defined as
$$ H^{m,n}_{p,q}(z) = H^{m,n}_{p,q}
\left[z   \left\vert
{\hfill (a_j ,\alpha _j )_{j=1,\dots,,p} \hfill \atop
 \hfill (b_k,\beta _k)_{k=1, \dots,q} \hfill}
 \right.  \right]=
\rec{2\pi i}\,
\int_L {\cal {H}}^{m,n}_{p,q} (s) \, z^s \, ds\,.\eqno(A.15)$$
We do not pursue furthermore in our survey:
we refer the interested reader to the treatises
on Fox $H$ functions
by  Mathai and Saxena \cite{MathaiSaxena H},
Srivastava, Gupta and Goyal \cite{Srivastava H} and
references therein.

\end{document}